\documentclass[11pt]{article}
\usepackage{enumerate}
\usepackage{amssymb,a4wide,latexsym,makeidx,epsfig,fleqn}
\usepackage{amsthm}
\usepackage{amsmath}
\usepackage{enumerate}
\usepackage{graphicx}
\usepackage{float}
\usepackage{color}
\usepackage{indentfirst}
\usepackage{tikz}
\usepackage{graphicx}
\usepackage{pdflscape}
\usepackage{siunitx} 
\usepackage{diagbox} 
\usepackage{rotating}
\usepackage{booktabs} 
\usepackage[colorlinks, linkcolor=blue, anchorcolor=blue, citecolor=blue]{hyperref}
\newtheorem{theorem}{Theorem}[section]

\newtheorem{lemma}[theorem]{Lemma}
\newtheorem{proposition}[theorem]{Proposition}

\newtheorem{problem}[theorem]{Problem}

\begin{document}
	\textwidth 150mm \textheight 225mm
	\title{Spectral radius and size conditions for 1-tough balanced bipartite graphs to contain a $k$-factor with $k\ge 3$
		\thanks{Supported by the National Natural Science Foundation of China (No. 12271439)}}
	
	\author{{Jinxuan Ren$^{a,b}$, Jiayao Song$^{a,b}$, Ligong Wang$^{a,b,}$\thanks{Corresponding author.}  }\\
		{\small $^a$School of Mathematics and Statistics, Northwestern
			Polytechnical University,}\\ {\small  Xi'an, Shaanxi 710129,
			P.R. China.}\\
		{\small $^b$ Xi'an-Budapest Joint Research Center for Combinatorics, Northwestern
			Polytechnical University,}\\
		{\small Xi'an, Shaanxi 710129,
			P.R. China. }\\
		{\small E-mail: jxrenmath@163.com,  Jiayaosongmath@163.com, lgwangmath@163.com} }
	\date{}
	\maketitle	
	\begin{center}
		\begin{minipage}{120mm}
			\vskip 0.3cm
			\begin{center}
				{\small {\bf Abstract}}
			\end{center}
			{\small   The largest eigenvalue of the adjacency matrix $A(G)$ of a graph $G$ is called the spectral radius $\rho(G)$ of $G$. If $S$ is a vertex cut of a connected graph $G$, then $c(G-S)$ denotes the number of connected components of $G-S$. The toughness $t(G)$ of a graph $G$ is defined as $t(G) = \min\left\{\frac{|S|}{c(G-S)}: S \subset V(G), \text{ and } c(G-S) > 1\right\}$.  A graph \( G \) is \( t \)-tough if \( |S| \ge t c(G - S) \) for every proper subset \( S \subset V(G) \)  satisfying \( c(G - S) > 1 \). Combining the toughness and spectral conditions, we provide spectral radius and size conditions for the existence of $k$-factors in 1-tough balanced bipartite graphs. This settles a problem posed by Chen, Fan and Lin in [Discrete Appl. Math. {\bf 355} (2024), 223--231].
				
				 }
		\end{minipage}
		
	\end{center}
	
	%%%%%%%%%%%%%%%%%%%%%%%%%%%%%%%%%%%%%%%%%%%%%%%%%%%%%%%%%%%%% Part 1: Introduction %%%%%%%%%%%%%%%%%%%%%%%%%%%%%%%%%%%%%%%%%%%%%%%%%%%%%%%%%%%%%%%%%%%%%%%%%%
	\section{Introduction }
	\label{sec:ch6-introduction}
	% Para1: basic concept
		This paper considers only undirected simple graphs. For a graph $G$ of order $n$, $V(G)$ denotes the vertex set and $E(G)$ denotes the edge set of $G$. We write $e(G)=|E(G)|$ for the size of $G$, i.e., the number of edges in $G$. The set of neighbors of a vertex $v\in V(G)$ in $G$ is denoted by $N_G(v)$, and $d_G(v)=|N_G(v)|$ is the degree of $v$ in $G$. The adjacency matrix of $G$ is the $n\times n$ matrix $A(G)=(a_{uv})_{n\times n}$, where $a_{uv}=1$ if $u$ and $v$ are adjacent, and $a_{uv}=0$ otherwise. For $i=1,2,\dots,n$, let $\lambda_i(G)$ be the $i$-th largest eigenvalue of $A(G)$, and write $\lambda_i=\lambda_i(G)$ for brevity. In particular, the largest eigenvalue $\lambda_1(G)$ of $A(G)$, also denoted by $\rho(G)$, is called the spectral radius of $G$. A graph is bipartite if its vertex set can be partitioned into two subsets $X$ and $Y$ such that every edge has one endpoint in $X$ and the other in $Y$. We denote a bipartite graph $G$ with bipartition \((X,Y)\) by \(G[X,Y]\). If $G[X,Y]$ is simple and every vertex in $X$ is joined to every vertex in $Y$, then $G$ is called a complete bipartite graph, denoted $K_{m,n}$, where $|X|=m$ and $|Y|=n$. A balanced bipartite graph of order $n=2m$ is denoted by $G[X,Y]$ with $|X|=|Y|=m$. A graph $H$ is called a subgraph of $G$ if $V(H)\subseteq V(G)$ and $E(H)\subseteq E(G)$. In particular, a graph $H$ is called a spanning subgraph of $G$ if $V(H)=V(G)$ and $E(H)\subseteq E(G)$. Let $g,f$ be two integer-valued functions defined on $V(G)$. A spanning subgraph $H$ of $G$ is known as a $(g,f)$-factor if $g(v)\leq d_H(v)\leq f(v)$ for every vertex $v\in V(G)$. When $g(x)=a$ and $f(x)=b$ for all $x\in V(G)$, such a $(g,f)$-factor is termed an $[a,b]$-factor. Furthermore, it is called an $f$-factor if $g(x)=f(x)$ for all $x\in V(G)$. In the special case where $a=b=k$, the $[a,b]$-factor reduces to a $k$-factor.
	
		Research on graph factors originated from the work of Danish mathematician Petersen \cite{Petersen} in 1891. He verified that any graph with all vertices of even degree can be partitioned into a collection of edge-disjoint $2$-factors. Meanwhile, Petersen \cite{Petersen} showed that every $2$-connected cubic graph contains a $1$-factor. These two landmark theorems laid the foundation for factor theory. In light of the characteristics of induced subgraphs, factor-related problems are generally classified into two categories: degree-constrained factors and component factors. A degree-constrained factor is a spanning subgraph of a graph such that the degree of each vertex satisfies prescribed lower and upper bounds. Typical examples of such factors include odd factors, even factors, $k$-factors and $[a,b]$-factors (see \cite{FanD,Kano}). A component factor is a spanning subgraph $H$ of a graph $G$ such that each connected component of $H$ is isomorphic to some member of a prescribed family of connected graphs. A number of concrete component factors have been extensively investigated. For example, an $\mathcal{S}_k$-factor is a component factor whose components are stars chosen from $\{K_{1,1},K_{1,2},\dots,K_{1,k}\}$, and a $P_{\ge k}$-factor is another component factor whose components are paths of order at least $k$. The above two categories are not mutually exclusive. As an illustration, a $1$-factor of a graph can be described as a spanning subgraph in which every vertex has degree exactly $1$, or equivalently, a spanning subgraph every connected component of which is the path $P_2$. Factor problems often intersect with other topics in graph theory, including the Hamiltonian cycle problem and the subgraph packing problem. Readers interested in further knowledge about graph factors may consult the monographs \cite{Aki,Yu}. 
	
		In recent decades, the study of conditions guaranteeing the existence of graph factors has attracted extensive attention, and a variety of approaches based on degree conditions and spectral conditions have been developed. Among these, characterizing factor existence from the perspective of graph spectra has emerged as a particularly active research direction, establishing a profound connection between factor theory and spectral graph theory. The investigation of eigenvalue conditions for factor existence can be traced back to the work of Brouwer and Haemers \cite{Brou2} in 2005, who gave a description of regular graphs containing a perfect matching in terms of the third largest adjacency eigenvalue. Subsequently, O \cite{SO} extended this result to general (not necessarily regular) graphs by providing a spectral condition to guarantee the existence of a perfect matching. For $k$-factors in general graphs, Lu \cite{Lu} obtained such a condition based on the third largest adjacency eigenvalue by employing Tutte's $k$-factor theorem. By incorporating the edge connectivity condition, Gu \cite{GuX1} studied eigenvalue conditions for $t$-edge-connected regular graphs containing $k$-factors. Tang and Zhang \cite{Tang} established spectral radius conditions for the existence of $k$-factors in graphs. Fan and Lin~\cite{FanDD} gave spectral radius conditions for the existence of $k$-factors in balanced bipartite graphs. More recently, Fan et al.~\cite{FanD} presented a comprehensive survey on eigenvalues and factors, summarizing relevant results on factors from the perspective of adjacency eigenvalues. Subsequent works have further extended these spectral results to various factor types, including $[a,b]$-factors \cite{Wei,FDD}, component factors \cite{ChenH,Jia,Haoyf}, even factors \cite{JiaC,LRF,Jhc}, even factors \cite{XZZ,ZSC,FA,HYL} and star factors \cite{Zhou,RFY,WSF}, among others.
	
		Chvátal \cite{Chv} first proposed the concept of graph toughness in 1973 . He defined the toughness $t(G)$ of a non-complete graph $G$ as $t(G)=\min\left\{\frac{|S|}{c(G-S)}\right\},$ where the minimum is taken over all proper subsets $S \subset V(G)$, and $c(G-S)$ denotes the number of connected components of $G-S$. This definition is closely related to many important graph properties. For example, it laid a foundation for the investigation of the relationship between graph toughness and the existence of factors \cite{Bian,Chen,Eno1,Eno2,Eno3}. It can be observed that the toughness $t(G)$ of any bipartite graph $G[X,Y]$ is at most $1$. The reason is that deleting all vertices in partition $X$ from $G$ (assuming $|X|\ge |Y|$) leaves $Y$ as an independent set. A non-complete bipartite graph $G$ is $t$-tough if $|S| \ge t\, c(G - S)$ for every proper subset $S \subset V(G)$ with $c(G - S) > 1$. Owing to the special structural properties of bipartite graphs, Bian \cite{Bian} further introduced the notion of bipartite toughness in 2006. The bipartite toughness of a non-complete balanced bipartite graph $G$ is defined as $t^B(G)=\min\left\{\frac{|S|}{c(G-S)}\right\},$ where the minimum is taken over all proper subsets $S\subset X$ (or $S\subset Y$) such that $G-S$ is disconnected and $c(G-S)>1$. Liu et al. \cite{Liu1} proved that a balanced bipartite graph of order $2n$ contains a $k$-factor if $n\geq4k-4$ and $t^B(G)>\frac{2(k-1)(n-1)}{kn+1}$, establishing a connection between bipartite toughness and the existence of $k$-factors. For studies concerning the relation between toughness and eigenvalues, we refer the reader to \cite{Alon,Brou,GuX}. Naturally, Chen, Fan and Lin considered the following problem in \cite{ChenY}.
	
		\begin{problem}(\cite{ChenY})\label{Pro1}
			What is the condition for 1-tough balanced bipartite graphs contains a $k$-factor for $k \ge 2$?
		\end{problem}
	
		For any two subsets $A,B\subseteq V(G)$, $E_G(A,B)$ denotes the set of edges in $G$ with one endpoint in $A$ and the other in $B$, and we write $e_G(A,B)=|E_G(A,B)|$. Given bipartite graphs $G_1[A_1,B_1]$ and $G_2[A_2,B_2]$, the notations $G_1[A_1,B_1]\:\nabla_1\:  G_2[A_2,B_2]$ and $G_1[A_1,B_1]\:\nabla_2\: G_2[A_2,B_2]$ represent graphs constructed from the union $G_1[A_1,B_1]\:\cup\: G_2[A_2,B_2]$ by adding all possible edges between $A_1$ and $B_2$, between $A_2$ and $B_1$, respectively.
	
\begin{figure}[htbp]
			\centering
			\tikzset{
				vertex/.style={circle, fill=black, inner sep=0pt, minimum size=7pt},
				edot/.style={circle, fill=black, inner sep=0pt, minimum size=2.2pt}
			}
			\begin{tikzpicture}[line width=1pt,scale=0.9]
				% 上方左侧椭圆 + k-1个顶点 u_1,...,u_{k-1}
				\draw (-3.2,3.8) ellipse (2.3cm and 0.68cm);
				\node[font=\large,above] at (-3.2,4.5) {$k-1$};
				\node[vertex,label=above:$u_1$] (u1) at (-4.8,3.6) {};
				\node[vertex,label=above:$u_{k-1}$] (uk1) at (-1.6,3.6) {};
				\node[edot] at (-3.9,3.6){};
				\node[edot] at (-3.3,3.6){};
				\node[edot] at (-2.7,3.6){};
			
				% 上方右侧椭圆 + n/2 -k+1个顶点 u_k,...,u_n
				\draw (3.2,3.8) ellipse (2.3cm and 0.68cm);
				\node[font=\large,above] at (3.2,4.5) {$\frac n2 -k+1$};
				\node[vertex,label=above:$u_k$] (uk) at (1.6,3.6) {};
				\node[vertex,label=above:$u_{n/2}$] (un) at (4.65,3.6) {};
				\node[edot] at (2.5,3.6){};
				\node[edot] at (3.1,3.6){};
				\node[edot] at (3.7,3.6){};
				
				% 下方左侧小圆 单个顶点 v_1
				\draw (-4.8,0.86) ellipse (0.68cm and 0.68cm);
				\node[font=\large,below] at (-4.8,0) {$1$};
				\node[vertex,label=below:$v_1$] (v1) at (-4.8,1) {};
			
				% 下方右侧椭圆 n/2 -1个顶点 v_2,...,v_n
				\draw (2.65,0.86) ellipse (3cm and 0.68cm);
				\node[font=\large,below] at (3.0,0) {$\frac n2 -1$};
				\node[vertex,label=below:$v_2$] (v2) at (0.45,1) {};
				\node[vertex,label=below:$v_{n/2}$] (vn) at (4.65,1.0) {};
				\node[edot] at (1.4,1){};
				\node[edot] at (2.5,1){};
				\node[edot] at (3.6,1){};
			
				% 连边：左上部{u1~uk1} 全连 v1、v2~vn；右上部{uk~un}全连v2~vn
				% u1连线
				\draw (u1)--(v1);
				\draw (u1)--(v2);
				\draw (u1)--(vn);
				% uk-1连线
				\draw (uk1)--(v1);
				\draw (uk1)--(v2);
				\draw (uk1)--(vn);
				% uk连线
				\draw (uk)--(v2);
				\draw (uk)--(vn);
				% un连线
				\draw (un)--(v2);
				\draw (un)--(vn);
			\end{tikzpicture}
			\caption{$G^0\big(k-1,\frac{n}{2}-k+1;\,1,\frac{n}{2}-1\big)$}
			\label{fig1}
		\end{figure}
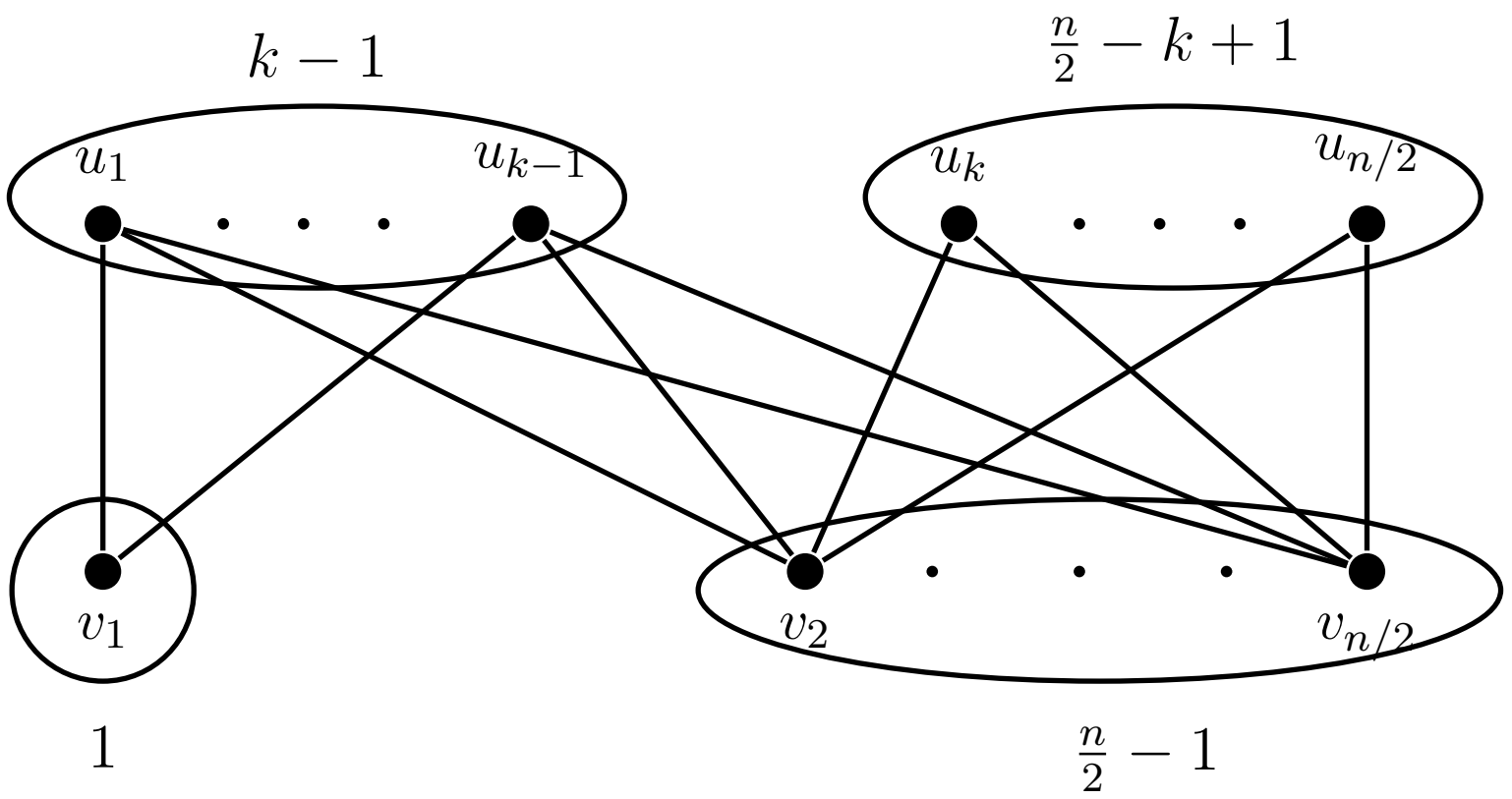

		Let $O_{a,b}=G[A,B]$ be a bipartite graph satisfying $|A|=a$, $|B|=b$ and $e_{O_{a,b}}(A,B)=0$. Denote by $O_{a,a}^{+a}$ the bipartite graph obtained by adding $a$ mutually disjoint edges to $O_{a,a}$. Define $G^{|X_1|}(|X_2|,|X_3|;|Y_2|,|Y_3|)=\big(O_{|X_1|,|Y_1|}^{+|X_1|}\nabla_2 K_{|X_2|,|Y_2|}\big)\nabla_1 K_{|X_3|,|Y_3|}$, where $|X_1|=|Y_1|$. Define $G^0(|X_2|,|X_3|;|Y_2|,|Y_3|)= K_{|X_2|,|Y_2|}\nabla_1 K_{|X_3|,|Y_3|}$, where $|X_1|=|Y_1|=0$. In particular, the graph $G^0\big(k-1,\frac{n}{2}-k+1;\,1,\frac{n}{2}-1\big)$ is illustrated in Figure \ref{fig1}.
	
		Chen, Fan and Lin \cite{ChenY} combined spectral radius with toughness and established a spectral radius condition for the existence of $2$-factors in $1$-tough balanced bipartite graphs, and proposed the following open problem.
	
		\begin{problem}(\cite{ChenY})\label{Pro2}
			Let $G$ be a $1$-tough balanced bipartite graph of order $n$. If
			\[
			\rho(G) \geq \rho\left(G^3\left(1, \frac{n}{2} - 4; 0, \frac{n}{2} - 3\right)\right),
			\]
			then $G$ contains a Hamiltonian cycle, unless $G \cong G^3\left(1, \frac{n}{2} - 4; 0, \frac{n}{2} - 3\right)$
		
		\end{problem}
	
		Ai and Zhang \cite{Ai} solved Problem \ref{Pro2}, and Chen, Fan and Lin \cite{ChenY} also proposed Problem \ref{Pro3}.
	
		\begin{problem}(\cite{ChenY})\label{Pro3}
			What is the spectral radius condition for $1$-tough balanced bipartite graphs with a $k$-factor for $k \geq 3$?
		
		\end{problem}
		In this paper, we will give a spectral radius condition and a size condition related to Problem \ref{Pro3} (see Theorems \ref{thm1} and \ref{thm2}).

		\begin{theorem}\label{thm1}
			Suppose that $G$ is a connected $1$-tough balanced bipartite graph of order $n$. If
			\[
			\rho(G)\ge \rho\left(G^0\left(k-1,\frac{n}{2}-k+1;\,1,\frac{n}{2}-1\right)\right), \: where \:  3\le k\le\frac{n}{2},
			\]
			then $G$ contains a $k$-factor, unless $G\cong G^0\left(k-1,\frac{n}{2}-k+1;\,1,\frac{n}{2}-1\right)$.
			
		\end{theorem}
	
		Theorem \ref{thm1} provides a sufficient spectral radius condition for $1$-tough balanced bipartite graphs containing a $k$-factor for $k\ge 3$. This condition is not necessary, because we can construct a $1$-tough balanced bipartite graph containing a $k$-factor whose spectral radius is smaller than $\rho\left(G^0\left(k-1,\frac{n}{2}-k+1;\,1,\frac{n}{2}-1\right)\right)$.
	
		\begin{proposition}\label{pos1}
			Let $k\ge 3$ be an integer. Then there exists a $1$-tough balanced bipartite graph $G^k\!\left(k-1,\frac{n}{2}-2k+1;0,\frac{n}{2}-k\right)$ containing a $k$-factor and
			\[
			\rho\!\left(G^k\!\left(k-1,\frac{n}{2}-2k+1;0,\frac{n}{2}-k\right)\right)<\rho\!\left(G^0\!\left(k-1,\frac{n}{2}-k+1;1,\frac{n}{2}-1\right)\right).
			\]
		
		\end{proposition}
	
		Proposition \ref{pos1} illustrates that it is impossible to characterize the necessary and sufficient spectral radius condition for $1$-tough balanced bipartite graphs  to contain a $k$-factor merely by the spectral radius comparison of a single extremal graph.
	
		We also consider Problem \ref{Pro3} from the size version, and provide the size condition to guarantee the existence of a $k$-factor ($k\ge 3$) in $1$-tough balanced bipartite graphs.
	
		\begin{theorem}\label{thm2}
			Suppose that $G$ is a connected $1$-tough balanced bipartite graph of order $n$. If
			\[
			e(G) \ge \frac{n^2}{4} - \frac{n}{2} + k - 1, \: where \:  3\le k\le\frac{n}{2},
			\]
			then $G$ contains a $k$-factor, unless $G\cong G^0\left(k-1,\frac{n}{2}-k+1;\,1,\frac{n}{2}-1\right)$.
		
		\end{theorem}
	
	%%%%%%%%%%%%%%%%%%%%%%%%%%%%%%%%%%%%%%%%%%%%%%%%%%%%%%%%%%%%% Part 2: Preliminaries %%%%%%%%%%%%%%%%%%%%%%%%%%%%%%%%%%%%%%%%%%%%%%%%%%%%%%%%%%%%%%%%%%%%%%%%%
	\section{Preliminaries}\label{sec:Preliminaries}
		In this section, we introduce some notations and several lammas that will be used throughout the paper.
	
		Let $M$ be an $n\times n$ real‑valued matrix, and let $X=\{1,2,\dots,n\}$. Suppose $\Pi: X=X_1\cup X_2\cup\cdots\cup X_k$ constitutes a partition of $X$. Under this partition, $M$ can be decomposed conformally into block form
		\[
		M=
			\begin{pmatrix}
				M_{1,1} & M_{1,2} & \cdots & M_{1,k}\\
				M_{2,1} & M_{2,2} & \cdots & M_{2,k}\\
				\vdots  & \vdots  & \ddots & \vdots\\
				M_{k,1} & M_{k,2} & \cdots & M_{k,k}
			\end{pmatrix}.
		\]
		The quotient matrix associated with the partition $\Pi$ is a $k\times k$ matrix $B_{\Pi}=(b_{i,j})_{i,j=1}^{k}$, where each entry $b_{i,j}$ equals the mean of all row sums taken from the sub‑block $M_{i,j}$. A partition $\Pi$ is said to be equitable whenever every block $M_{i,j}$ of $M$ possesses identical row sum $b_{i,j}$. Correspondingly, we refer to $B_{\Pi}$ as an equitable quotient matrix provided that $\Pi$ forms an equitable partition for $M$.
	
		\begin{lemma}[\textup{\cite{Brou1,Godsil}}]\label{lem1}
			Let $M$ be a real symmetric matrix, and let $\lambda_1(M)$ be the largest eigenvalue of $M$. If $B_{\Pi}$ is an equitable quotient matrix of $M$, then the eigenvalues of $B_{\Pi}$ are also eigenvalues of $M$. Furthermore, if $M$ is nonnegative and irreducible, then $\lambda_1(M) = \lambda_1(B_{\Pi})$.
		\end{lemma}
	
		By the well-known Perron--Frobenius theorem \cite{Godsil}, we can easily deduce the following result.
	
		\begin{lemma}[\textup{\cite{ChenY}}]\label{lem2}
			If $H$ is a spanning subgraph of a connected graph $G$, then
			\[
			\rho(H)\leq \rho(G),
			\]
			with equality if and only if $H\cong G$.
		\end{lemma}
	
		\begin{lemma}[\textup{\cite{Liu1}}]\label{lem3}
			Let $G=G[X,Y]$ be a balanced bipartite graph, and $k\geq 1$ be an integer. Then $G$ has a $k$-factor if and only if for any $S\subseteq X$ and $T\subseteq Y$,
			\[
			k|S|+\sum_{v\in T}d_G(v)-k|T|-e_G(S,T)\geq 0.
			\]
		\end{lemma}
	
		\begin{lemma}\label{lem4}
			If $\dfrac{n}{2}\ge k$, then the complete balanced bipartite graph $K_{\frac{n}{2},\frac{n}{2}}$ contains a $k$-factor.
		\end{lemma}
	
		\begin{proof}
			Let $G=K_{\frac{n}{2},\frac{n}{2}}$ with bipartition $(X,Y)$, where $|X|=|Y|=\frac{n}{2}$.
			By Lemma \ref{lem3}, We have that the balanced bipartite graphs $G$ contains a $k$-factor if and only if for any subsets $S\subseteq X$ and $T\subseteq Y$,
			\[
			k|S|+\sum_{v\in T}d_G(v)-k|T|-e_G(S,T)\ge 0.
			\]
			Denote $|S|=s$, $|T|=t$. Since each vertex of  $G=K_{\frac{n}{2},\frac{n}{2}}$ has degree $d_G(v)=\frac{n}{2}$ and $e_G(S,T)=st$, it suffices to verify that for any subsets $S\subseteq X$ and $T\subseteq Y$,
			\[
			ks+\frac{n}{2}t-kt-st=s(k-t)+t(\frac{n}{2}-k)\geq 0,
			\]
			Let $\Phi(s,t)=ks+\frac{n}{2}t-kt-st=s(k-t)+t\left(\frac{n}{2}-k\right).$
			where $0\le s,t\le \frac{n}{2}$ and $\frac{n}{2}\ge k$.
		
			By symmetry of the two partitions of a complete bipartite graph, it suffices to treat the cases for $t$. Those for $s$ follow by symmetry.
		
			{\textbf{Case 1.}}  $t=0$.
		
			In this case, $\Phi(s,0)=ks\ge0$.
		
			{\textbf{Case 2.}}  $0<t\le k$.
		
			In this case, we have $k-t\ge0$ and $\frac{n}{2}-k\ge0$. Since $s\ge0$, both terms are non‑negative, so $\Phi(s,t)=s(k-t)+t\left(\frac{n}{2}-k\right)\ge0$.
		
			{\textbf{Case 3.}}  $t>k$. 
		
			In this case, we can consider $\Phi(s,t)$ as a linear function of $s$:
			\[
			f(s)=(k-t)s+t(\frac{n}{2}-k).
			\]
			Since $k-t<0$, $f(s)$ is decreasing in $s$. The minimum attains at $s=\frac{n}{2}$:
			\[
			f(\frac{n}{2})=(k-t)\frac{n}{2}+t(\frac{n}{2}-k)=k(\frac{n}{2}-t)\ge 0,
			\]
			Thus $f(s)\ge0$ for all $0\le s\le \frac{n}{2}$.
			
			All cases yield $\Phi(s,t)\ge0$. Then $K_{\frac{n}{2},\frac{n}{2}}$ contains a $k$-factor.
		\end{proof}

	%	\begin{corollary}\label{lem7}
		%	Let $G$ be a connected graph of even order $n\geq 4$. If $G-\{u,v\}$ has a perfect $k$-matching for any 
		%$	u,v \in V(G)$, then $G$ is $\text{GBC}_k$.
		%	\end{corollary}

	%%%%%%%%%%%%%%%%%%%%%%%%%%%%%%%%%%%%%%%%%%%%%%%%%%%%%%%%%%%%% Part 3: Proofs %%%%%%%%%%%%%%%%%%%%%%%%%%%%%%%%%%%%%%%%%%%%%%%%%%%%%%%%%%%%%%%%%%%%%%%%%
	\section{ Proofs of Main Results }
	
	Firstly, we prove Theorem \ref{thm1}
	
	\noindent
		{\textbf{Proof of Theorem \ref{thm1}}.
		Let $G^*$ be the connected balanced bipartite graph constructed by adding any one edge to $G^0\big(k-1,\frac{n}{2}-k+1;\,1,\frac{n}{2}-1\big)$. From the structure of the graph in Figure \ref{fig1}, it follows that any added edge must be chosen among $v_1u_k, v_1u_{k+1},\dots, v_1u_{n/2}$. By Lemma \ref{lem2}, it follows that $\rho(G^*) > \rho\big(G^0(k-1,\frac{n}{2}-k+1;\,1,\frac{n}{2}-1)\big)$. 
		
		Let $G$ be a connected $1$-tough balanced bipartite graph. We have the following three claims.
		
		{\parindent=0pt
			{\textbf{Claim 1.}}\label{claim1}
			$\delta(G)\ge 2$.
		}
		
		\begin{proof}
			
			Otherwise, suppose that $\delta(G)\le 1$.
			
			{\textbf{Case 1.}} $\delta(G)=0$.
			
			In this case, $G$ contains an isolated vertex, which contradicts the connectivity of a tough graph, a contradiction.
			
			{\textbf{Case 2.}} $\delta(G)=1$.
			
			In this case, there exists a vertex $u\in V(G)$ with $d_G(u)=1$. Let $v$ be the unique neighbor of $u$.
			Take the vertex cut $S=\{v\}$. After removing $S$, $u$ becomes an isolated component, so $c(G-S)\ge 2$. Then
			$\frac{|S|}{c(G-S)}=\frac{1}{c(G-S)}\le \frac12<1$, a contradiction.

			By Cases 1 and Case 2, Claim 1 holds.
		\end{proof}
		
		{\parindent=0pt
			{\textbf{Claim 2.}}\label{claim2}
			Both the graph $G^0\!\left(k-1,\frac{n}{2}-k+1;1,\frac{n}{2}-1\right)$ and the above $G^*$ are  $1$-tough balanced bipartite graph. 
		}

		\begin{proof}
			Write $G^0\!\left(k-1,\frac{n}{2}-k+1;1,\frac{n}{2}-1\right) = G^0$. By construction, $|X|=|X_2|+|X_3|=\frac{n}{2}$ and $|Y|=|Y_2|+|Y_3|=\frac{n}{2}$, so $G^0$ is balanced bipartite.
			Write $Y_2=\{v_1\}$. The vertex $v_1$ is adjacent only to the vertices in $X_2$, and there are no edges between $v_1$ and $X_3$. Every vertex in $Y_3$ is adjacent to all vertices of $X$. Let $S$ be an arbitrary vertex cut of $G^0$, set $S_X=S\cap X$, $S_Y=S\cap Y$, and $c=c(G^0-S)$. It suffices to show $|S|\ge c$.
			
			We split the discussion into two cases depending on whether all vertices of $Y_3$ are removed by $S$. The set $Y_3$ acts as a connecting hub. If at least one vertex of $Y_3$ remains, the remaining vertices in $X_2$ and $X_3$ are connected. If $Y_3$ is completely removed, the vertices in $X_3$ lose all their neighbors and become isolated vertices.
			
			\textbf{Case 1.} $Y_3\setminus S_Y\neq \emptyset$.
			
			In this case, all surviving vertices in $X_2\setminus S_X$ and $X_3\setminus S_X$ are connected via vertices in $Y_3\setminus S_Y$, forming a large component. The only possible extra component is the one containing $v_1$.
			
			\textbf{Subcase 1.1.} If $v_1\in S_Y$, then $c\le 1$. Hence $|S|\ge c$.
			
			\textbf{Subcase 1.2.} If $v_1\notin S_Y$ and $X_2\setminus S_X\neq\emptyset$, then $v_1$ belongs to the large component. Hence $c\le 1$.
			
			\textbf{Subcase 1.3.} If $v_1\notin S_Y$ and $X_2\subseteq S_X$, then $v_1$ becomes an isolated vertex and $c=2$. Since $k\ge 3$, $|S_X|\ge |X_2|=k-1\ge 2$, which yields $|S|\ge 2=c$.
			
			\textbf{Case 2.} $Y_3\subseteq S_Y$.
			
			In this case, $Y\setminus S_Y\subseteq \{v_1\}$, and each vertex in $X_3\setminus S_X$ is isolated. Let $p=|X_3\setminus S_X|$. Then $c\le p+1\le \frac{n}{2}-k+2\le \frac{n}{2}-1$.
			Since $|S_Y|\ge |Y_3|=\frac{n}{2}-1$, we get $|S|\ge \frac{n}{2}-1\ge p+1\ge c$.
			
			Thus $|S|\ge c$ holds for every vertex cut $S$, so $G^0$ is $1$-tough.
			
			The only missing edges of $G^0$ are between $v_1$ and vertices in $X_3$. Therefore, the added edge for $G^*$ must be $v_1x_0$ for some $x_0\in X_3$.
			For any vertex cut $S$ of $G^*$, adding edges cannot increase the number of connected components, so $c(G^*-S)\le c(G^0-S)$. Combining the $1$-tough property of $G^0$,
			\[
			\frac{|S|}{c(G^*-S)}\ge \frac{|S|}{c(G^0-S)}\ge 1.
			\]
			Hence $G^*$ is also $1$-tough.
		\end{proof}
		
		{\parindent=0pt
			{\textbf{Claim 3.}}\label{claim3}
			$G^0\!\left(k-1,\frac{n}{2}-k+1;1,\frac{n}{2}-1\right)$ contains no $k$-factor.
		}
		
		\begin{proof}
			Since the degree of vertex $v_1$ in $G^0\!\left(k-1,\frac{n}{2}-k+1;1,\frac{n}{2}-1\right)$ is $k-1$, this graph cannot contains a $k$-factor.
		\end{proof}
		
		Since only undirected simple graphs are considered in this paper, it follows from the above Claims 1,2 and 3 that any edge added to $G^0\big(k-1,\frac{n}{2}-k+1;1,\frac{n}{2}-1\big)$ must be incident to the vertex $v_1$. Therefore, if $\rho(G)>\rho\big(G^0\big(k-1,\frac{n}{2}-k+1;1,\frac{n}{2}-1\big)\big)$, then $\rho(G)\ge \rho(G^*)$. It suffices to verify that $G^*$ contains a $k$-factor. 
		
		{\parindent=0pt
			{\textbf{Case 1. $\frac{n}{2}=k\ge 3$}}
		}
		
		In this case, since $G^*$ is obtained by adding an  edge to $G^0$, adding one edge to $G^0$ yields the complete bipartite graph $K_{k,k}$. It follows from Lemma \ref{lem4} that $G^*$ contains a $k$-factor.
		
		{\parindent=0pt
			{\textbf{Case 2. $\frac{n}{2}>k\ge 3$}}
		}
		
		Consider the complete bipartite graph $K_{\frac{n}{2}-1,\frac{n}{2}-1}$ with bipartition $X=\{u_2,u_3,\dots,u_{\frac{n}{2}}\}$ and $Y=\{v_2,v_3,\dots,v_{\frac{n}{2}}\}$. Since $\frac{n}{2}>k$, it follows that $\frac{n}{2}-1\ge k$. By Lemma \ref{lem4}, the complete bipartite graph $K_{\frac{n}{2}-1,\frac{n}{2}-1}$ must contain a $k$-factor. We arbitrarily select one such $k$-factor and denote the corresponding spanning subgraph by $H_k$. Then $H_k$ is a $k$-regular bipartite graph satisfying
		$d_{H_k}(u_2)=d_{H_k}(u_3)=\cdots=d_{H_k}(u_{\frac{n}{2}})=k,\:\:\:d_{H_k}(v_2)=d_{H_k}(v_3)=\cdots=d_{H_k}(v_{\frac{n}{2}})=k$. In $H_k$, the $k-1$ vertices $u_2,u_3,\dots,u_k$ are incident to a total of $k(k-1)$ edges. Now we introduce a new vertex $u_1$, and redistribute these $k(k-1)$ edges evenly among the $k$ vertices $u_1,u_2,\dots,u_k$. This construction yields a new bipartite graph $H'$, whose bipartition is given by $X=\{u_1,u_2,\dots,u_k,u_{k+1},\dots,u_{\frac{n}{2}}\},\:\:\: Y=\{v_2,v_3,\dots,v_{\frac{n}{2}}\}$.
		
		The vertex degrees in $H'$ satisfy
		\begin{align*}
			d_{H'}(u_1)=d_{H'}(u_2)=\cdots=d_{H'}(u_k)&=k-1,\\
			d_{H'}(u_{k+1})=d_{H'}(u_{k+2})=\cdots=d_{H'}(u_{\frac{n}{2}})&=k,\\
			d_{H'}(v_2)=d_{H'}(v_3)=\cdots=d_{H'}(v_{\frac{n}{2}})&=k.
		\end{align*}
		We further add a vertex $v_1$ to $H'$, and add the edges $u_1v_1,u_2v_1,\dots,u_kv_1$. The resulting graph $H''$ is a $k$-regular bipartite graph with $|X|=|Y|=\frac{n}{2}$.
		
		It can be observed that $H''$ is isomorphic to a $k$-factor of $G^*$, where $G^*$ denotes the graph obtained by adding an arbitrary edge to $G^0\!\left(k-1,\frac{n}{2}-k+1;1,\frac{n}{2}-1\right)$. Consequently, $G^*$ contains a $k$-factor isomorphic to $H''$.
		
		Combining the above arguments, $G^*$ contains a $k$-factor whenever $\frac{n}{2}\ge k\ge 3$. This completes the proof of Theorem \ref{thm1}.\nobreak\hfill$\square$
		
		Secondly, we prove Proposition \ref{pos1}}.
		
	\noindent
		{\textbf{Proof of Proposition \ref{pos1}}. We consider the graph $G^k\!\left(k-1,\frac{n}{2}-2k+1;0,\frac{n}{2}-k\right)$ and the graph $G^0\big(k-1,\frac{n}{2}-k+1;\,1,\frac{n}{2}-1\big)$. We assume that the orders of the two graphs are sufficiently large.
			
			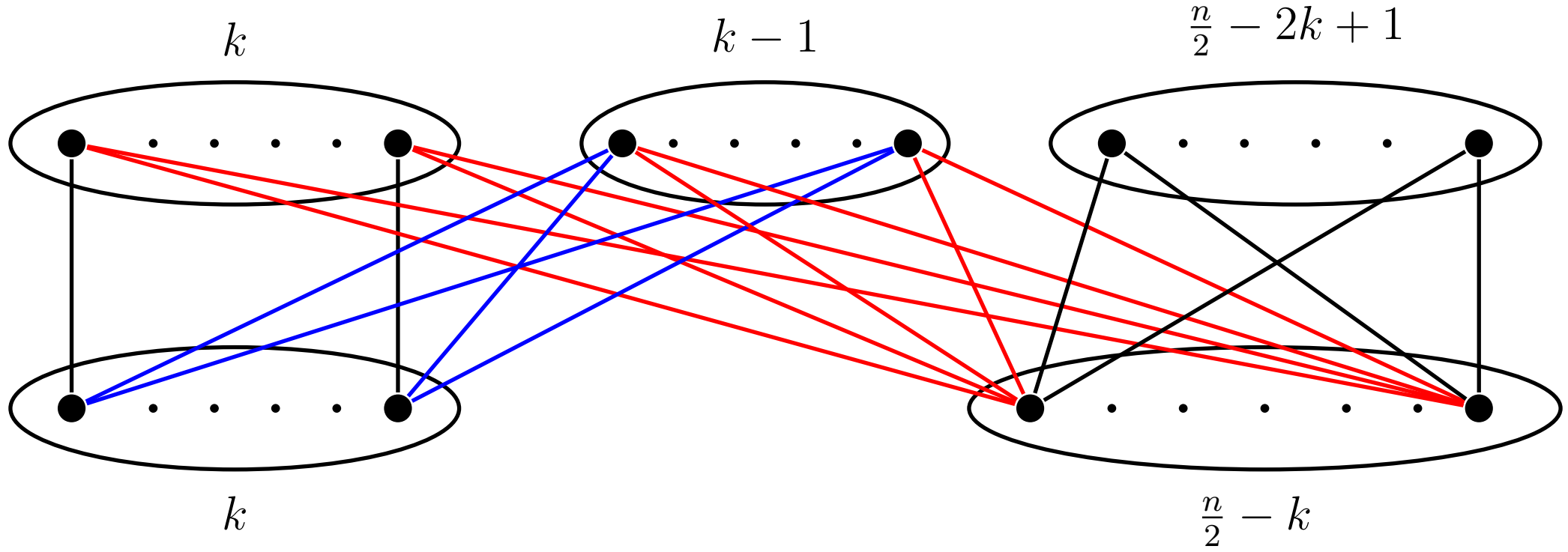
\begin{figure}[htbp]
				\centering
				\tikzset{
					vertex/.style={circle, fill=black, inner sep=0pt, minimum size=7pt},
					edot/.style={circle, fill=black, inner sep=0pt, minimum size=2.2pt}
				}
				\begin{tikzpicture}[line width=1pt,scale=0.9]
					% ========== 上方三个顶点集 ==========
					% 上左：k 个顶点
					\draw (-5.2,3.6) ellipse (2.2cm and 0.6cm);
					\node[font=\large,above] at (-5.2,4.3) {$k$};
					\node[vertex] (uL1) at (-6.8,3.6) {};
					\node[vertex] (uLk) at (-3.6,3.6) {};
					\node[edot] at (-6.0,3.6){};
					\node[edot] at (-5.4,3.6){};
					\node[edot] at (-4.8,3.6){};
					\node[edot] at (-4.2,3.6){};
					
					% 上中：k-1 个顶点
					\draw (0,3.6) ellipse (1.8cm and 0.6cm);
					\node[font=\large,above] at (0,4.3) {$k-1$};
					\node[vertex] (uM1) at (-1.4,3.6) {};
					\node[vertex] (uMkm1) at (1.4,3.6) {};
					\node[edot] at (-0.9,3.6){};
					\node[edot] at (-0.3,3.6){};
					\node[edot] at (0.3,3.6){};
					\node[edot] at (0.9,3.6){};
					
					% 上右：\frac n2 -2k+1 个顶点
					\draw (5.2,3.6) ellipse (2.4cm and 0.6cm);
					\node[font=\large,above] at (5.2,4.3) {$\frac n2 -2k+1$};
					\node[vertex] (uR1) at (3.4,3.6) {};
					\node[vertex] (uRt) at (7.0,3.6) {};
					\node[edot] at (4.1,3.6){};
					\node[edot] at (4.7,3.6){};
					\node[edot] at (5.4,3.6){};
					\node[edot] at (6.1,3.6){};
					
					% ========== 下方两个顶点集 ==========
					% 下左：k 个顶点
					\draw (-5.2,1) ellipse (2.2cm and 0.6cm);
					\node[font=\large,below] at (-5.2,0.3) {$k$};
					\node[vertex] (vL1) at (-6.8,1.0) {};
					\node[vertex] (vLk) at (-3.6,1.0) {};
					\node[edot] at (-6.0,1.0){};
					\node[edot] at (-5.4,1.0){};
					\node[edot] at (-4.8,1.0){};
					\node[edot] at (-4.2,1.0){};
					
					% 下右：\frac n2 -k 个顶点
					\draw (4.9,1) ellipse (2.9cm and 0.6cm);
					\node[font=\large,below] at (4.8,0.3) {$\frac n2 -k$};
					\node[vertex] (vR1) at (2.6,1.0) {};
					\node[vertex] (vRs) at (7.0,1.0) {};
					\node[edot] at (3.4,1.0){};
					\node[edot] at (4.1,1.0){};
					\node[edot] at (4.9,1.0){};
					\node[edot] at (5.7,1.0){};
					\node[edot] at (6.4,1.0){};
					
					% ========== 绘制边，模仿原图三色示意：黑色、红色、蓝色 ==========
					% 上左集合(uL1,uLk) 连下左(vL1,vLk) 黑色竖边
					\draw (uL1)--(vL1);
					\draw (uLk)--(vLk);
					
					% 上左集合 红色边连向下右集合两端
					\draw[red] (uL1)--(vR1);
					\draw[red] (uL1)--(vRs);
					\draw[red] (uLk)--(vR1);
					\draw[red] (uLk)--(vRs);
					
					% 上中集合(uM1,uMkm1)：蓝色边连下左，红色边连下右
					\draw[blue] (uM1)--(vL1);
					\draw[blue] (uM1)--(vLk);
					\draw[blue] (uMkm1)--(vL1);
					\draw[blue] (uMkm1)--(vLk);
					
					\draw[red] (uM1)--(vR1);
					\draw[red] (uM1)--(vRs);
					\draw[red] (uMkm1)--(vR1);
					\draw[red] (uMkm1)--(vRs);
					
					% 上右集合(uR1,uRt) 黑色交叉边只连下右集合
					\draw (uR1)--(vR1);
					\draw (uR1)--(vRs);
					\draw (uRt)--(vR1);
					\draw (uRt)--(vRs);
					
				\end{tikzpicture}
				\caption{$G^k\!\left(k-1,\frac{n}{2}-2k+1;0,\frac{n}{2}-k\right)$}
				\label{fig2}
			\end{figure}
			
			We have the following three claims
			
			{\parindent=0pt
				{\textbf{Claim 4.}}\label{claim4}
				$G^k\!\left(k-1,\frac{n}{2}-2k+1;0,\frac{n}{2}-k\right)$ is a $1$-tough balanced bipartite graph.
			}
			
			\begin{proof}
				From the structure illustrated in Fig. \ref{fig2}, it is easy to verify that there exists no subset $S\subset G^k\!\left(k-1,\frac{n}{2}-2k+1;0,\frac{n}{2}-k\right)$ such that $|S|<c\big(G^k\!\left(k-1,\frac{n}{2}-2k+1;0,\frac{n}{2}-k\right)-S\big)$. Hence, the toughness of $G^k\!\left(k-1,\frac{n}{2}-2k+1;0,\frac{n}{2}-k\right)$ is $1$.
			\end{proof}

			{\parindent=0pt
				{\textbf{Claim 5.}}\label{claim5}
				$\rho\!\left(G^k\!\left(k-1,\frac{n}{2}-2k+1;0,\frac{n}{2}-k\right)\right)<\rho\!\left(G^0\!\left(k-1,\frac{n}{2}-k+1;1,\frac{n}{2}-1\right)\right)$.
			}
			
			\begin{proof}
				Let $G^k=G^k\!\left(k-1,\frac{n}{2}-2k+1;0,\frac{n}{2}-k\right)$ and $G^0=G^0\!\left(k-1,\frac{n}{2}-k+1;1,\frac{n}{2}-1\right)$. Then $A(G^k)$ has the equitable quotient matrix $B_{\Pi_1}$, and $A(G^0)$ has the equitable quotient matrix $C_{\Pi_2}$.
				\[
				B_{\Pi_1}=
				\begin{pmatrix}
					0 & 0 & 0 & 1 & \frac{n}{2}-k \\
					0 & 0 & 0 & k & \frac{n}{2}-k \\
					0 & 0 & 0 & 0 & \frac{n}{2}-k \\
					1 & k-1 & 0 & 0 & 0 \\
					k & k-1 & \frac{n}{2}-2k+1 & 0 & 0
				\end{pmatrix},\quad
				C_{\Pi_2}=
				\begin{pmatrix}
					0 & 0 & 1 & \frac{n}{2}-1 \\
					0 & 0 & 0 & \frac{n}{2}-1 \\
					k-1 & 0 & 0 & 0 \\
					k-1 & \frac{n}{2}-k+1 & 0 & 0
				\end{pmatrix}
				\]
				The characteristic polynomial of $B_{\Pi_1}$ is
				\[
				\begin{aligned}
					\mu_1(B_{\Pi_1},x)
					&= x^5 - ( \frac14 n^2 -\frac12 kn + k^2 - k + 1)x^3 \\
					&\quad - (-k^4 + k^3 n - \frac12 k^2 n - \frac14 k^2 n^2 + \frac14 kn^2 - \frac14 n^2)x, \\[6pt]
				\end{aligned}
				\]
				and the characteristic polynomial of $C_{\Pi_2}$ is
				
				\[
				\begin{aligned}
					\mu_2(C_{\Pi_2},x)
					&= x^4 + (\frac12 n - \frac14 n^2 - k + 1)x^2
					- \frac{k^2 n}{2} + k^2 + \frac{kn}{2} \\
					&\quad -2k + 1 + \frac{kn^2}{4} - \frac{n^2}{4}.
				\end{aligned}
				\]
				
				When $n$ is sufficiently large and $x>\sqrt{\frac{n}{2}\left(\frac{n}{2}-1\right)}$, we have $\mu_1(B_{\Pi_1},x)-\mu_2(C_{\Pi_2},x)>0$, and hence $\mu_1(B_{\Pi_1},x)>\mu_2(C_{\Pi_2},x)$ for $x>\sqrt{\frac{n}{2}\left(\frac{n}{2}-1\right)}$. Since $G^0$ contains $K_{\frac{n}{2},\frac{n}{2}-1}$ as a proper subgraph, $\rho(G^0) > \rho(K_{\frac{n}{2},\frac{n}{2}-1}) = \sqrt{\frac{n}{2}\left(\frac{n}{2}-1\right)}$ by lemma \ref{lem2}. Thus, $\lambda_1(C_{\Pi_2}) > \lambda_1(B_{\Pi_1})$. According to Lemma \ref{lem1}, we have  $\rho(G^0) > \rho(G^k)$.
			\end{proof}
			
			{\parindent=0pt
				{\textbf{Claim 6.}}\label{claim6}
				$G^k\!\left(k-1,\frac{n}{2}-2k+1;0,\frac{n}{2}-k\right)$ contains a $k$-factor.
			}
			
			\begin{proof}
				Let $H_0[X_0,Y_0]=O_{|X_1|,|Y_1|}^{+|X_1|}\nabla_2 K_{|X_2|,|Y_2|}$ be a bipartite graph, where $X_1=\{u_1,u_2,...,u_k\}$, $Y_1=\{v_1,v_2,...,v_k\}$, $X_2=\{u_{k+1},u_{k+2},...,u_{2k-1}\}$, $Y_2=\emptyset$, $X_0 = X_1\:\cup\: X_2$, $Y_0= Y_1\:\cup\: Y_2$. Let $H_4[X_4,Y_4]= K_{|X_4|,|Y_4|}$ be a bipartite graph, where $X_3=\{u_{2k},u_{2k+1}...,u_{\frac{n}{2}}\}$, $X_4=X_1\:\cup\:X_3$, $Y_4=\{v_{k+1},v_{k+2},...,v_{\frac{n}{2}}\}$. 
				
				We observe that $H_4[X_4,Y_4] - \{u_1\} \cong K_{\frac{n}{2}-k,\frac{n}{2}-k}$.
				For sufficiently large $n$, by Lemma \ref{lem4}, the graph $H_4[X_4,Y_4] - \{u_1\}$ contains a $k$-factor.
				Take one such $k$-factor and denote it by $H_5[X_5,Y_5]$, or simply $H_5$.
				
				By the definition of $k$-factor, the vertex degrees in $H_5$ satisfy
				\[
				d_{H_5}(u_2)=d_{H_5}(u_3)=\cdots=d_{H_5}(u_k)=d_{H_5}(u_{2k})=\cdots=d_{H_5}(u_{\frac{n}{2}})=k,
				\]
				\[
				d_{H_5}(v_{k+1})=d_{H_5}(v_{k+2})=\cdots=d_{H_5}(v_{\frac{n}{2}})=k.
				\]
				
				Therefore,  $e(H_5) = k\left(\frac{n}{2}-k\right)$.
				We have the algebraic identity $k\left(\frac{n}{2}-k\right)=k(k-1)+\left(\frac{n}{2}-2k+1\right)k$. With this edge-count identity, we may locally rearrange the edges of $H_5$ to obtain a new bipartite graph $H_6[X_6,Y_6]$ satisfying $X_6 = \{u_1,u_2,\dots,u_k,u_{2k},\dots,u_{\frac{n}{2}}\},\:\:Y_6 = \{v_{k+1},v_{k+2},\dots,v_{\frac{n}{2}}\}$, along with the degree conditions below
				\[
				\begin{aligned}
					d_{H_6}(u_1) &= d_{H_6}(u_2) = \cdots = d_{H_6}(u_k) = k-1,\\
					d_{H_6}(u_{2k}) &= d_{H_6}(u_{2k+1}) = \cdots = d_{H_6}(u_{\frac{n}{2}}) = k,\\
					d_{H_6}(v_{k+1}) &= d_{H_6}(v_{k+2}) = \cdots = d_{H_6}(v_{\frac{n}{2}}) = k.
				\end{aligned}
				\]
				
				Set $H'[X,Y]=H'[X_ 0\cup X_6,Y_0 \cup Y_6]$. The vertex set of $H'[X,Y]$ is $X = X_0 \cup X_6=\{u_1,...,u_k,u_{k+1},...,u_{2k-1},u_{2k},...,u_{\frac{n}{2}}\}$ and $Y = Y_0 \cup Y_6=\{v_1,...,v_k,v_{k+1},...,v_{\frac{n}{2}}\}$, and the edge set of $H'[X,Y]$ is $E(H_0[X_0,Y_0] \cup H_6[X_6,Y_6])$. Then $H'[X,Y]$ is a $k$-regular balanced bipartite graph with $|X|=\frac{n}{2},\quad |Y|=\frac{n}{2}$.
				We compute the vertex degrees as follows
				\[
				\begin{aligned}
					&d_{H'}(u_1)=d_{H'}(u_2)=\cdots=d_{H'}(u_k)=d_{H_0}(u_1)+d_{H_6}(u_1)=d_{H_0}(u_2)+d_{H_6}(u_2)\\
					&=\cdots=d_{H_0}(u_k)+d_{H_6}(u_k)=1+k-1=k,
				\end{aligned}
				\]
				\[
				d_{H'}(u_{k+1})=d_{H'}(u_{k+2})=\cdots=d_{H'}(u_{2k-1})=d_{H_0}(u_{k+1})=d_{H_0}(u_{k+2})=\cdots=d_{H_0}(u_{2k-1})=k,
				\]
				\[
				d_{H'}(u_{2k})=d_{H'}(u_{2k+1})=\cdots=d_{H'}(u_{\frac{n}{2}})=d_{H_6}(u_{2k})=d_{H_6}(u_{2k+1})=\cdots=d_{H_6}(u_{\frac{n}{2}})=k,
				\]
				\[
				d_{H'}(v_1)=d_{H'}(v_2)=\cdots=d_{H'}(v_k)=d_{H_0}(v_1)=d_{H_0}(v_2)=\cdots=d_{H_0}(v_k)=k,
				\]
				\[
				d_{H'}(v_{k+1})=d_{H'}(v_{k+2})=\cdots=d_{H'}(v_{\frac{n}{2}})=d_{H_6}(v_{k+1})=d_{H_6}(v_{k+2})=\cdots=d_{H_6}(v_{\frac{n}{2}})=k.
				\]
				
				From the construction above, $H'$ is a spanning subgraph of $G^k\left(k-1,\frac{n}{2}-2k+1;0,\frac{n}{2}-k\right)$.
				Consequently, $G^k\left(k-1,\frac{n}{2}-2k+1;0,\frac{n}{2}-k\right)$ contains a $k$-factor.
			\end{proof}
			
			From the results of Claims 4, 5 and 6, we consider the two $1$-tough balanced bipartite graphs $G^k\left(k-1,\frac{n}{2}-2k+1;0,\frac{n}{2}-k\right)$
			and $G^0\left(k-1,\frac{n}{2}-k+1;1,\frac{n}{2}-1\right).$
			Observe that their spectral radii satisfy $
			\rho\left(G^k\right)<\rho\left(G^0\right),$
			where $G^0\left(k-1,\frac{n}{2}-k+1;1,\frac{n}{2}-1\right)$ contains no $k$-factor, whereas
			$G^k\left(k-1,\frac{n}{2}-2k+1;0,\frac{n}{2}-k\right)$ contains a $k$-factor.
			
			This completes the proof.\nobreak\hfill$\square$
			
			Finally, we prove Theorem \ref{thm2}}.
			
			\noindent
				{\textbf{Proof of Theorem \ref{thm2}}}.
				Since $e\left(G^0\!\left(k-1,\frac{n}{2}-k+1;1,\frac{n}{2}-1\right)\right) = \frac{n^2}{4}-\frac{n}{2}+k-1$ and $G^0\!\left(k-1,\frac{n}{2}-k+1;1,\frac{n}{2}-1\right)$ contains no $k$-factor. Therefore, it suffices to prove that if $e(G)\ge \frac{n^2}{4}-\frac{n}{2}+k$, then every $1$-tough balanced bipartite graph $G$ of order $n$ contains a $k$-factor.
				
				We proceed by contradiction. Suppose that $G=G[X,Y]$ has no $k$-factor. By Lemma \ref{lem3}, there exist subsets $S\subseteq X$ and $T\subseteq Y$ such that
				\begin{equation}
					k|S|+\sum_{v\in T}d_G(v)-k|T|-e_G(S,T)<0. \tag{1}
				\end{equation}
				Since $G$ is bipartite, we have $\sum_{v\in T}d_G(v)=e_G(S,T)+e_G(X\setminus S,T)$. Substituting this into (1), we obtain $k|S|+e_G(X\setminus S,T)-k|T|<0$,
				that is,
				\begin{equation}
					e_G(X\setminus S,T)<k(|T|-|S|). \tag{2}
				\end{equation}
				
				Let $|S|=s$ and $|T|=t$. Since $e_G(X\setminus S,T)\ge 0$, Inequality (2) implies $t>s$. Moreover, $e_G(X\setminus S,T)$ is an integer, so
				\begin{equation}
					e_G(X\setminus S,T)\le k(t-s)-1. \tag{3}
				\end{equation}
				
				We now estimate the missing sizes of $G$ relative to the complete bipartite graph $K_{n/2,n/2}$. By establishing a lower bound for these missing sizes, we obtain an upper bound on the size of $G$, which yields a contradiction.
				
				Note that $e_{K_{n/2,n/2}}\left(X\setminus S,T\right) = (\frac{n}{2}-s)t$. For the graph $G$, only $e_G(X\setminus S,T)$ sizes exist between these two subsets. Hence the missing sizes for this bipartite pair are at least $(\frac{n}{2}-s) t - e_G(X\setminus S,T)$. Using Inequality (3), we derive a lower bound for the missing sizes
				\begin{equation}
					(\frac{n}{2}-s)t - e_G(X\setminus S,T) \ge (\frac{n}{2}-s)t - k(t-s)+1. \tag{4}
				\end{equation}
				
				We have the following claim.
				
				\noindent
					{\textbf{Claim 7.}}\label{claim7}
					Suppose that $0\le s \le \frac{n}{2}-1$, $1\le t-s \le \frac{n}{2}-s$, and $k\le \frac{n}{2}$. Then
					\begin{equation}
						(\frac{n}{2}-s)t - k(t-s) \ge \frac{n}{2}-k. \tag{5}
					\end{equation}
					\begin{proof}
						Let $r = t-s.$ Then $1\le r \le \frac{n}{2}-s$. Substitute $t = s+r$ into the inequality (5)
						\begin{align}
							(\frac{n}{2}-s)t - k(t-s)
							&= (\frac{n}{2}-s)(s+r) - kr \notag \\
							&= s(\frac{n}{2}-s) + r(\frac{n}{2}-s-k)\notag\\
							&\geq \frac{n}{2}-k \tag{6}
						\end{align}
						
						We consider two cases depending on the value of $s$.
						
						\medskip
						\textbf{Case 1.} $s \le \frac{n}{2}-k$. 
						
						In this case, we have $\frac{n}{2}-s-k \ge 0$. Since $r\ge 1$, we have $r(\frac{n}{2}-s-k) \ge \frac{n}{2}-s-k$. We obtain the following inequality
						\begin{align*}
							s(\frac{n}{2}-s)+r(\frac{n}{2}-s-k)
							&\ge s(\frac{n}{2}-s) + (\frac{n}{2}-s-k) \\
							&= (\frac{n}{2}-s)(s+1)-k\\
							&\ge \frac{n}{2}-k.\:\:(\mbox{Since } 0 \le s \le \frac{n}{2}-1)
						\end{align*}
						
						Thus, in this case, Inequality (6) holds.
						
						\medskip
						\textbf{Case 2.} $s > \frac{n}{2}-k$. 
						In this case, we have $\frac{n}{2}-s-k < 0$. Since $r \le \frac{n}{2}-s$, we have $r\left(\frac{n}{2}-s-k\right) \ge \left(\frac{n}{2}-s\right)(\frac{n}{2}-s-k)$.
						We obtain the following inequality:
						\begin{align*}
							s(\frac{n}{2}-s)+r(\frac{n}{2}-s-k)
							&\ge s(\frac{n}{2}-s)+(\frac{n}{2}-s)(\frac{n}{2}-s-k) \\
							&= (\frac{n}{2}-s)(s+\frac{n}{2}-s-k) \\
							&= (\frac{n}{2}-s)(\frac{n}{2}-k)\\
							&\ge \frac{n}{2}-k.\:\:(\mbox{Since }s \le \frac{n}{2}-1)
						\end{align*}
						Thus, in this case, Inequality (6) holds.
						
						Combining Cases 1 and 2, Claim 7 holds. 
					\end{proof}
					
					By Claim 7, $(\frac{n}{2}-s)t - k(t-s) \ge \frac{n}{2}-k$. 
					Thus, inequality (4) yields $(\frac{n}{2}-s)t - e_G(X\setminus S,T) \ge \frac{n}{2}-k+1$. This means that $G$ is missing at least $\frac{n}{2}-k+1$ edges compared with the complete bipartite graph $K_{n/2,n/2}$, so $e(G) \le \frac{n^2}{4}-(\frac{n}{2}-k+1) = \frac{n^2}{4}-\frac{n}{2}+k-1$. This contradicts $e(G) \ge \frac{n^2}{4}-\frac{n}{2}+k$. 
					
					This completes the proof.\nobreak\hfill$\square$ 
					
					%Other methods for generating tables with continuous borders
					%	\begin{table}[htbp]
						%	
						%	
						%		%\footnotesize
						%			\centering
						%		\caption{Adjacency spectral radius of $G_s$ and $K_\frac{n}{2}\vee \frac{n}{2}K_1$}\label{example}
						%		%%\begin{tabular}{|c|*{19}{|c}|}
						%		\begin{tabular}{|c|*{4}{|c}|}
							%			\toprule
							%			$\rho(K_\frac{n}{2}\vee \frac{n}{2}K_1)$ &\diagbox{$n$}{$G_s$} & $\rho(G_1)$&$\rho(G_2)$&$\rho(G_3)$ \\
							%			\midrule
							%			2.56 &$n=4$& 2.17 & & \\
							%			4.16 &$n=6$& 4.05 & 3.62 &\\
							%			5.77 &$n=8$& 6.02 & 5.27 & 5.18\\
							%		
							%		
							%			\bottomrule
							%	\end{tabular}
						%	\label{tab:spectral_data}
						%	
						%\end{table}
						%	
						%	
						%	
						%	\begin{table}[htbp]
							%		\centering
							%		\caption{Adjacency spectral radius of $G_s$ and $K_\frac{n}{2}\vee \frac{n}{2}K_1$}
							%		\label{example}
							%		\begin{tabular}{ccccc}
								%			\toprule
								%			$\rho(K_\frac{n}{2}\vee \frac{n}{2}K_1)$ & \diagbox{$n$}{$G_s$} & $\rho(G_1)$ & $\rho(G_2)$ & $\rho(G_3)$ \\
								%			\midrule
								%			2.56 & $n=4$ & 2.17 & & \\
								%			4.16 & $n=6$ & 4.05 & 3.62 & \\
								%			5.77 & $n=8$ & 6.02 & 5.27 & 5.18 \\
								%			\bottomrule
								%		\end{tabular}
							%		\label{tab:spectral_data}
							%	\end{table}
						
	\section{ Concluding remarks}
						
		In this paper, we present a sufficient spectral radius condition for the existence of a $k$-factor in $1$-tough balanced bipartite graphs when $k\ge 3$.
		Moreover, we show that it is impossible to characterize the necessary and sufficient spectral radius condition for $1$-tough balanced bipartite graphs to contain a $k$-factor merely by the spectral radius comparison of a single extremal graph. 
		We remark that the extremal graph $G^0\left(k-1,\frac{n}{2}-k+1;1,\frac{n}{2}-1\right)$
		constructed in this paper fails to be $1$-tough for $k=2$.
		The sufficient condition for the existence of a $k$-factor in $1$-tough balanced bipartite graphs with $k=2$ has already been established in \cite{ChenY}.
						
		Furthermore, we also establish a size condition for the existence of a $k$-factor in 1-tough balanced bipartite graphs. We observe that the extremal graphs employed for the proposed spectral radius condition and size condition coincide.

	\section*{Declaration of competing interest}
						
		The authors declare that they have no known competing financial interests or personal relationships that could have appeared to influence the work reported in this paper.
						
	\section*{Data availability}
						
		No data was used for the research described in the article.

\end{document}